\documentclass[11pt]{article}
\newcommand{\mf }{{\mathfrak M}}

\newcommand{\vx}{\underline{x}}

\newcommand{\va}{\underline{a}}
\newcommand{\vb}{\underline{b}}

\newcommand{\vd}{\underline{\delta}}

\usepackage{diagrams}
\usepackage{amscd}
\usepackage{amsmath}
\usepackage{amssymb}
\usepackage{graphicx}
\usepackage{mathptmx}
\usepackage{amsfonts}
\newtheorem{prop}{Proposition}
\newtheorem{defn}[prop]{Definition}
\newtheorem{thm}[prop]{Theorem}

\begin{document}

\title{ Theory of characteristics for  first order partial 
differential equations}
\author{Anders Kock}
\date{}
\maketitle

%\bigskip

\small
\begin{flushright}{\em Dabei sehen wir von unendlich kleinen 
Gr\"{o}ssen h\"{o}here Ordnung ab.}\\
Lie \cite{BT} p.\ 523\end{flushright}
\normalsize
\section*{Introduction}
The present note makes no claim of originality; it is a ``conspectus'' of some of the classical theory 
of characteristics for 1st order PDEs, as expounded geometrically by 
Lie and elaborated by Klein. These authors use extensively a synthetic geometric 
language, but ultimately describe notions rigourously only by presenting 
them in analytic terms. Our approach describes the notions (like 
``united position'' (``vereinigte Lage'') and ``characteristic'') 
rigour\-ously in pure  
synthetic coordinate free terms, and introduces coordinates only at a later 
point, when it comes to proving some of the relations between the 
notions introduced.

So we are not  claiming that describing the notions synthetically is 
an effective tool for {\em proving}; usually, coordinates are better 
suited for this. The virtue of the synthetic descriptions are, as also 
appears from the work of Monge, Lie, Klein, \ldots , that it gives a geometric 
language to speak about geometric entities, and in particular, 
make them coordinate free from the outset.

The particular version of synthetic language that we use is that of 
Synthetic Differential Geometry, as in \cite{SDG}, say,  and notably  
as in \cite{SGM}, where the main synthetic relation is the first and 
second order neighbour relation, as first considered in French algebraic 
geometry in the 1950s. We denote these relations by the symbol $\sim 
_{1}$ (or just $\sim$) and $\sim _{2}$, respectively. They are 
reflexive symmetric relations on the set of points of a manifold. The set of $k$th order 
neighbours of a point $x$ in a manifold $M$ is denoted $\mf _{k}(x)$, 
i.e.\ $\mf _{k}(x) = \{y\in M \mid y\sim _{k}x\}$. 
One has that $x\sim _{1}y \sim _{1}z$ implies $x\sim _{2}z$.
The axiomatics used 
for these neighbourhoods is essentially the ``Kock-Lawvere'' (KL) axiom 
scheme, which we shall quote when needed. The basic manifold is the 
number line $R$; here $x\sim _{k}y$ iff $(y-x)^{k+1}=0$. In $R^{n}$, 
the set $\mf _{1}(0)$ is also denoted $D(n)$, and $\mf _{k}(0) $ is 
denoted $D_{k}(n)$.

\section{Surface elements and calottes}\label{SEC}
Let $M$ be a 3-dimensional manifold.

A {\em surface element} at $x$ is a set $P\subset M$ of the form $\mf _{1}(x) 
\cap F$, where $F\subset M$ is a surface (2-dimensional submanifold) containing $x$. Similarly a 
{\em calotte} at $x$ is a set $K\subset M$ 
of the form $\mf _{2}(x)\cap F$ where $F\subset M$ is a surface containing $x$.
(The notion of calotte is from \cite{HG} p.\ 281.)
If $K$ is a calotte at $x$, it is clear that $\mf _{1}(x)\cap K$ is a 
surface element $P$ at $x$, called the {\em restriction} of $K$, and 
similarly, $K$ is an {\em extension} calotte of $P$, or a calotte 
{\em through} $P$.

It follows from Proposition \ref{basept} in the Appendix that the base point $x$ of 
a surface element $P$ can be reconstructed\footnote{It is possible 
that the synthetic ``combinatorics'' presented here makes sense in 
other contexts than Synthetic Differential Geometry; in that case, one might probably have to 
consider the base point $x$ of a surface element $P$ as part of the data 
of it.} from $P$ (viewed as a 
subset), and 
similarly, the restriction of a calotte $K$ (and hence also the base 
point of the calotte) can be reconstructed from $K$ (viewed as a 
subset). 

One could also use the terms ``1-jet (resp.\ 2-jet) of a surface''  
for surface elements, respectively calottes, in $M$.
Therefore, and for uniformity, we denote the manifold of surface elements, 
respectively the manifold of calottes, by the symbols $S_{1}(M)$, 
respectively $S_{2}(M)$. We have surjective submersions
$$S_{2}(M)\to S_{1}(M)\to M.$$
The dimensions of these manifolds are 8, 5, and 3, respectively, cf.\ 
Section \ref{CC}.
The manifold $S_{1}(M)$ may be described as the projectivization 
$P(T^{*}M)$ of 
the cotangent bundle $T^{*}(M)\to M$.

A calotte $K$ at $x$ defines a family of surface elements  
namely the family of sets $\mf _{1}(y)\cap K$ for $y$ ranging over 
$\mf _{1}(x)$. The surface elements $P'$ coming about from $K$ in this way 
are said to {\em belong to} $K$, or be {\em contained in} $K$. Note that the restriction of $K$ 
belongs to $K$; a surface element which belongs to $K$ is the 
restriction of $K$ iff its base point is $x$.

If $F\subset M$ is a surface, there is a map $F\to S_{1}(M)$, 
associating to $x\in F$ the surface element $\mf _{1}(x)\cap F$.

\section{The contact distribution $\approx$}\label{CD}

We consider a general 3-dimensional manifold $M$, and the 
corresponding 5-dimen\-sional manifold $S_{1}(M)$ of surface elements.

Being a manifold, $S_{1}(M)$ carries a (1st order) neighbour relation 
$\sim$. It carries a further structure, namely   a reflexive symmetric relation
$\approx $ refining 
$\sim$, and called ``united position'' (``vereinigte Lage'', in the 
terminology og Lie and Klein): if $P$ and $Q$ are  neigbour 
surface elements with base points $p$ and $q$, respectively, we say that
$$P \approx Q \mbox {\quad  if \quad }q\in P.$$
This is almost a literal translation of the definition in Lie 
\cite{BT} 
p.\ 523: ``a surface element is in united position with 
another one if the point of the latter lies in the plane of the 
former''. It is not immediate from the definition that $\approx$ is a 
{\em symmetric} relation, but this can be proved (see Section \ref{CC}) if we, 
in Lie's verbal rendering (\cite{BT} p.\ 523), {\em ``ignore infinitesimally small quantities 
of higher order''}. In our context, the ``ignored quantities'' are 
not only ignored, but are {\em equal} to 0, using $P\sim 
Q$, as the coordinate calculation below (beginning of Section 
\ref{CC}) will reveal.

 Let $F$ be a surface in $M$. Since the passage from 
points $x$ in $F$ to the corresponding surface elements $\mf 
_{1}(x)\cap F$ is a function, 
it follows from general principles that the surface elements of $F$ 
at $x$ and $y$ (both in $F$, and with $x\sim y$) are neighbours in $S_{1}(M)$. 
Furthermore, $y\in \mf _{1}(x) 
\cap F$; so $y$ belongs to the surface element of $F$ at $x$. 
Thus we see that if $F$ 
is a surface, the surface elements at neighbouring points of $F$ are 
in united position. This is the motivation for the notion.

It follows that  if $K$ is a calotte at $x$, and $P$ is a surface element 
belonging to the calotte, then $P\approx K_{1}$, where $K_{1}$ is 
the surface element obtained by restriction of $K$. Consider namely 
some 
surface $F$ such that $K = \mf _{2}(x)\cap F$, and apply the 
above reasoning to
$F$.

(In modern treatments, the structure ``united position'' is 
presented as subordinate to the canonical {\em contact manifold} structure which 
the cotangent bundle $T^{*}M$ carries -- a certain canonical 1-form.  
However, $P(T^{*}M)$ does not carry a canonical 1-form (only ``modulo 
a scalar factor''), and our 
description (i.e.\ Lie's) of $\approx$ is purely geometric.)

\section{First order PDEs}\label{PDE}

By a {\em first order PDE} on a 3-dimensional manifold $M$, one 
understands a 4-dimen\-sional submanifold $\Psi$ of the 5-dimensional 
manifold 
$S_{1}(M)$ of surface elements in $M$. The {\em 
solutions} of $\Psi$ are then the  surfaces $F$ in $M$ such that 
all surface elements of $F$ belong to $\Psi$.

This geometric formulation of the analytic notion of ``first order partial 
differential equation'' goes back to Monge, Lie, and other 19th 
century geometers, cf.\  classical texts like \cite{BT}, 
\cite{HG}, \cite{Ste},\ldots .

By a {\em solution calotte} of $\Psi$, we mean a calotte all of whose surface 
elements belong to $\Psi$.
A necessary condition that a calotte $K$ at $x$ is a solution calotte is 
of course that its restriction belongs to $\Psi$,
i.e.\ $\mf _{1}(x)\cap K \in \Psi$.
We ask the converse question: let $P\in \Psi$. How many solution 
calottes  through 
$P$ are there, i.e.\ how many solution calottes are there
 with restriction $P$? We shall  prove that the set of such 
calottes form a 1-dimensional manifold, see Section \ref{CC}.

\section{Characteristic neighbours}\label{CN}
A neighbour surface element $P'$ of $P$ belonging to {\em some} 
 solution calotte of $\Psi$ through $P$ 
(is in united position with $P$, and belongs to $\Psi$, but) may not 
belong to {\em all} solution calottes through $P$. 

We   ask: given $P \in \Psi$, how many neighbour surface 
elements $P'$ of $P$  have the 
property that they belong to {\em all} these $\infty ^{1}$ solution 
calottes through $P$?  We pose:

\begin{defn}Let $\Psi$ be a PDE, and let $P\sim P'$ be neighbour 
surface elements in $\Psi$. If $P'$ belongs to {\em all} solution calottes 
through $P$, we say that $P'$ is a {\em characteristic} neighbour of 
$P$, written $P\approx _{\Psi}P'$.\end{defn}
(Note that $P\approx _{\Psi}$ implies $P\approx P'$.) 
Thus, if $F$ is a solution surface of $\Psi$ and contains $P$, then 
$F$ will also contain $P'$.
In particular, if two solution surfaces $F_{1}$ and $F_{2}$ are tangent to 
each other at $x$, meaning that $x\in F_{1}\cap F_{2}$ and $\mf_{1}(x)\cap F_{1}=\mf_{1}(x)\cap 
F_{2}$ ($=P$, say), and if $P'$ is a 
characteristic neighbour of $P$, then $F_{1}$ and $F_{2}$ both contain $P'$, 
equivalently, the surfaces are tangent to each other at the base point of $P'$.

We shall prove that the characteristic neighbour relation $\approx 
_{\Psi}$  defines a 1-dimen\-sional 
distribution on the manifold $\Psi$, and hence can be integrated into 
curves. These curves are the classical ``characteristic stripes'' of the PDE 
$\Psi$.

Thus, if if two solution surfaces $F_{1}$ and $F_{2}$ are tangent to 
each other at $x$,  then they are tangent to each other along the 
characteristic stripe through $P$.

If $P$ and $P'$ are characteristic neighbours, and $x'$ is the base 
point of $P'$, then $P'$ can be reconstructed from $x'$ and $P$. For, 
take any solution calotte $K$ through $P$ (such calottes do exist 
- there are in fact $\infty ^{1}$ of them). Since $P\approx P'$, we 
have that 
$x'\in P \subseteq K$, and since $P'$ belongs to all such solution 
calottes by assumption,  $P'=\mf _{1}(x')\cap K$ (and this is independent 
of the choice of $K$).

A point $x'$ which appears as the base point of a characteristic 
neigbour $P'$ of $P$ may be called a  characteristic neighbour {\em 
point} of $P$ ``in the calotte sense''. There is (for $M=R^{3}$) another, older, notion of characteristic 
neighbour point of $P$, going back to Monge, Lagrange, \ldots , 
namely, it is a point $x'$ of $P$, on the line along which $P$ is tangent to the 
``Monge cone'' at $x$ (where $x$ is the base point of $P$). We shall describe 
these notions in synthetic form in Section \ref{MC}, and prove that $x'$ is 
a characteristic neighbour point of $P$ in the ``calotte'' sense iff it is 
so in the ``Monge'' sense. This we have been unable to  prove from 
the purely synthetic data, and we prove it by establishing the 
differential equations that analytically express the synthetic notions of 
``characteristic''.

\section{Monge cone}\label{MC}
In the classical treatment, the manifold $M$ is $R^{3}$, and the 
surface elements in $R^{3}$ are called {\em plane} elements, since a surface 
element at $\vx \in R^{3}$ may be given by a {\em plane} through 
$\vx$. The plane elements of a PDE $\Psi$ through a fixed point $\vx$ 
have an enveloping surface, which is a cone, called the Monge cone at 
$\vx$; 
each individual plane element $P \in \Psi$ through $\vx$ is 
tangent to the Monge cone at $\vx$ along a generator of the cone, and 
this generator $l\subseteq P$ is the {\em characteristic line} of the plane 
element. Paraphrasing, we then arrive at the provisional definition 
that $\vx '$ is a characteristic neighbour (in the ``Monge sense'')  of the plane element $P\in 
\Psi $ 
through $\vx$ if $\vx ' \in \mf _{1}(\vx )\cap l$.

However, as argued in \cite{END}, the relationship between 
enveloping surfaces and characteristics is that the characteristics 
are logically prior to the enveloping surface (which is made up of the 
characteristics). From this conception, it is therefore a detour to 
define the characteristic lines $l$ in terms of the Monge cones. In 
fact, we define directly the notion of characteristic neighbour (``in 
the Monge sense''), and applicable for any PDE $\Psi$ on a 3-dimensional manifold $M$.
(The set of characteristic neighbours of $P$, as $P$ ranges over those 
surface elements in $\Psi$ which have base point $x$, is then an 
infinitesimal version of the classical Monge cone at $x$.)

\begin{defn} Let $\Psi$ be a PDE on a 3-dimensional manifold $M$, and let $P\in \Psi$ with 
base point $x$. Then $x'\in P$ is a characteristic neighbour for $P$ 
(in the ``Monge'' sense) if for all $P'\sim P$ with $x\in P' \in \Psi$, we have $x'\in P'$.
\end{defn}
This may be seen as a rigourous formulation of the description of 
Lie, \cite{BT} p.\ 510:  ``--{\em \ldots so hat man im Punkte 
$(x,y,z)$ die Schnittlinie der Ebenen zweier solcher unmittelbar 
benachbarte Fl\"{a}chenelemente \ldots zu suchen \ldots''} (he is 
talking about two plane elements through $(x,y,z) $. So  
instead of intersecting $P$ with ``an 
immediate neighbour'' $P'$, we intersect it with {\em all} its 
``immediate'' (first order) neigbours; this is the key idea in the 
conception of \cite{END}.)

\medskip

From the synthetic considerations in the previous Section, it is clear what 
role characteristic neigbour points of $P$, in the calotte sense, have for 
solutions. The role of the charcteristic neighbour points in the 
Monge sense is not immediately clear, but  we  prove, by analytic 
means (cf.\ the end of Section \ref{CC})  
that the notions agree. So we get that  the  
characteristic neighbours, in the Monge sense,  have the same synthetic 
role for solutions of the PDE $\Psi$ as, more evidently, the 
characteristic neigbours in the calotte sense do:
\begin{thm}\label{monge}If $x'\in P \in \Psi$ is a characteristic 
neighbour for the surface element $P$, and $F_{1}$ and $F_{2}$ are 
solutions of $\Psi$ containing $P$, then $\mf _{1}(x')\cap F_{1}=\mf 
_{1}(x')\cap F_{2}$.
\end{thm}

\section{Coordinate calculations}\label{CC}

We consider the case where $M=R^{3}$. A function $f:R^{2}\to R$ gives 
rise to a surface $F$ in $R^{3}$, namely its graph. Not all surfaces in 
$R^{3}$ come about this way (they may contain vertical surface 
elements), but since our considerations are local, it suffices to 
consider such ``graph''-surfaces.

If $(x,y)\in R^{2}$ and $f:R^{2}\to R$ is a function with graph $F$, 
then the restriction of $f$ to 
$\mf _{1}(x,y)$ has for its graph a 
surface element of $F$ at $(x,y,f(x,y))$. By KL axiom, this 
restriction is determined by $f(x,y)$ and by the two partial 
derivatives $p=\partial f /\partial x (x,y)$, 
$q=\partial f /\partial y (x,y)$, and is therefore a synthetic 
rendering of the 1-jet of $f$ at $(x,y)$. Thus,
the surface element determines the 5-tuple
$(x,y,z,p,q)$.

 Similarly, the 2-jet of $f$ 
at $(x,y)$ is the restriction of $f$ to $\mf _{2}(x,y)$; its graph  
is a calotte at $(x,y)$,
and it determines the 8-tuple $(x,y,z,p,q,r,s,t)$ with $z,p,q$ 
as before, and with $r,s,t$ the second order partial derivatives of $f$ at 
$(x,y)$, $r= \partial ^{2}f /\partial x^{2}$, $s= \partial ^{2}f 
/\partial x\partial y$, and $t= \partial ^{2}f /\partial y^{2}$ 
(evaluated at $(x,y)$).

Conversely,  any 5-tuple $(x,y,z,p,q)$ defines a (non-vertical) surface 
element $P$ in $R^{3}$, namely the graph 
of the 1-jet at $(x,y)$ of 
the  affine function $f_{1}:R^{2}\to R$ given  by
\begin{equation}f_{1}(\xi ,\eta )= z + p(\xi -x)+q(\eta 
-y).\label{PQ}\end{equation} 
\begin{sloppypar}Similarly, any
8-tuple $(x,y,z,p,q,r,s,t)$ defines a calotte, namely 
 the graph of the 2-jet at $(x,y)$ of the quadratic 
function
$f_{2}: R^{2}\to R$ given by  
\begin{equation}f_{2}(\xi ,\eta )= z + p(\xi -x)+q(\eta -y) + 
\tfrac{1}{2}r(\xi -x)^{2} 
+s(\xi -x)(\eta -y) + \tfrac{1}{2}t(\eta 
-y)^{2}.\label{RST}\end{equation}\end{sloppypar}
(Note that the 1-jet of  the function $f_{2}$ at $(x,y)$ agrees with 
the 1-jet at $(x,y)$ of the function $f_{1}$, since on $\mf _{1}(x,y)$, the second order terms 
vanish.)

The points belonging to the surface element $(x,y,z,p,q)$ are the points of the form
$$(x+dx, y+dy, z+p\; dx + q\; dy)$$
with $(dx,dy)\in D(2)$.\footnote{We follow Klein and Lie in this 
notation for ``first order infinitesimal elements'', i.e.\ for elements 
in $D, D(2),\ldots$, but we want to emphasize that $dx, dy, \ldots $ are {\em not} differential 
forms (which behave contravariantly), but rather,  $dx$ and $dy$ are 
elements of $R$ (``numbers''), behaving in a certain sense {\em 
covariantly}; more precisely, the neighbour relation $\sim$ is 
preserved by any map between manifolds.} The base point of this surface element is 
$(x,y,z)$.

The points  belonging to the calotte $(x,y,z,p,q,r,s,t)$ are the 
points of the form
    $$(x+\delta x, y+\delta y, z+p\; \delta x + q\; \delta y + 
    \tfrac{1}{2}r(\delta x)^2 + s\; \delta x \delta y + \tfrac{1}{2} 
    t(\delta y)^2)$$
    with $(\delta x ,\delta y) \in D_{2}(2)$. 

\medskip

For two surface elements $P$ and $P'$ to be in united position, 
$P\approx P'$,  they must  
first of all be neighbours, $P\sim P'$, so they are of the form
$$(x,y,z,p,q)\mbox{\quad  and \quad }(x+dx, y+dy,z+dz,p+dp,q+dq)$$
respectively, with $(dx,dy,dz,dp,dq)\in D(5)$; and then 
\begin{equation}P\approx P' \mbox{\quad iff \quad}dz = p\; dx + q\; dy.
    \label{upo}\end{equation}

To prove symmetry of the relation $\approx$, we should  
from $dz = p\; dx + q\;dy$ deduce that
$$-dz = (p+dp)(-dx)+(q+dq)(-dy);$$
but this follows because $dp\cdot dx =0$ and $dq\cdot dy =0$ since 
$(dx, \ldots ,dq)\in D(5)$.  (Lie puts it this way, \cite{BT} p.\ 523: ``{\em here, 
we ignore infinitely small quantities of higher order}''; in our 
formalism, the ``higher order quantities'' to be ignored are 
$dp\cdot dx$ and $dq\cdot dy$; they are both 0.) 

\medskip

We consider a surface element $P=(x,y,z,p,q)$ and ask for the 
relation between on the one hand

$\bullet$ calottes $K=(x,y,z,p,q,r,s,t)$ 
extending $P$, and

$\bullet$ surface elements $P'= (x+dx,y+dy,z+p\; dx +q\; dy, 
p+dp, q+dq)$ in united position with $P$, on the other. (Here, 
$(dx,dy,dp,dq)\in D(4)$.)

\begin{prop}\label{PPrime}The surface element $P'$ belongs to $K$  
iff $(r,s,t)$ is a solution of a certain linear equation 
system (two equations in  three unknowns), namely
 the 
linear system with augmented matrix  
\begin{equation}\left[\begin{array}{ccc|c}
dx&dy&&dp \\
&dx&dy&dq \end{array}\right].
\label{surff2}\end{equation}
\end{prop} 
{\bf Proof.}Consider the function $f=f_{2}$ from (\ref{RST}), whose 
2-jet at $(x,y)$ has $K$ as graph. Its first partial derivatives at 
$(x+dx,y+dy)$ are $p+r\; dx +s\; dy$ and  $q+s\; dx + t\; dy$, 
respectively. For $¬P'$ to belong to $K$, these partial derivatives 
have to be $p+dp$ and $q+dq$, respectively.
 This equation expresses a relation between $(r,s,t)$ and 
$(dx,dy,dp,dq)$ on the other, which may be rewritten in matrix form as 
stated.

\medskip

\begin{sloppypar}Now we bring in a PDE $\Psi$, a 4-dimensional submanifold of the 
5-dimen\-sional manifold $S_{1}(R^{3})$ of surface elements in $R^{3}$. Our 
considerations are local, so we may assume that $\Psi$ is given as 
the zero set of a certain function $\psi :R^{5} \to R$, in other 
words,
$(x,y,z,p,q)\in \Psi$ iff $\psi (x,y,z,p,q)=0$. The graph of a 
function $f:R^{2}\to R$ is then a solution surface iff for all 
$(x,y)$
$$\psi (x,y,f(x,y), \tfrac{\partial f}{\partial x}(x,y),\tfrac{\partial 
f}{\partial y}(x,y))=0$$
which is a  partial differential equation of order 1 (and this is the 
justification for our more general use of the term ``PDE''). We assume that 
``$p$ and $q$ really occur in the function $\psi$'' (or ``$\psi$ is 
not free of $p$ and $q$''):  we 
assume that 
$\partial \psi /\partial p$ and $\partial \psi /\partial q$ do not 
vanish simultaneously at any point $(x,y,z,p,q)$ (more precisely: at 
any $(x,y)$, at least one of $\partial \psi /\partial p$ and $\partial 
\psi /\partial q$ is invertible).\end{sloppypar}

We proceed to describe the solution calottes for $\Psi$ in analytic 
terms. A necessary condition that a calotte $K=(x,y,z,p,q,r,s,t)$ is 
a solution calotte is of course that its restriction $(x,y,z,p,q)$ is 
in $\Psi$.

\begin{prop}\label{Psii}Assume $(x,y,z,p,q)\in \Psi$. Then 
the calotte $(x,y,z,p,q,r,s,t)$ is a solution calotte for 
$\Psi$ iff $(r,s,t)$ is a solution of the linear equation system (two 
equations in  three unknowns), with 
augmented matrix
\begin{equation}\left[\begin{array}{ccc|c}
\psi _{p}&\psi _{q}&&-\psi _{x} - p\cdot \psi _{z}\\
&\psi _{p}&\psi_{q} &-\psi _{y} -q\cdot \psi _{z} \end{array}\right].
\label{surff4}\end{equation}
where $\psi _{x}$ denotes $\frac{\partial \psi}{\partial 
x}(x,y,z,p,q)$, and similarly for $\psi _{y}, \psi _{z}, \psi 
_{p},\psi _{q}$.
 \end{prop}
{\bf Proof.} Let $f=f_{2}:R^{2}\to R$ be the quadratic function given 
by (\ref{RST}). The calotte in question is then a solution calotte 
iff for all $(dx,dy)\in D(2)$
\begin{equation*}\psi (x+dx,y+dy,z+p\; dx +q\; dy , \partial f/\partial x, \partial 
f/\partial y) =0\end{equation*}
where the partial derivatives are to be evaluated at $(x+dx,y+dy)$; 
these partial derivatives are $r\cdot dx + s\cdot dy$ and $s\cdot dx + t\cdot dy$, 
respectively, so $K$ is a solution calotte iff
\begin{equation}\psi (x+dx,y+dy,z+p\; dx +q\; dy ,p+r\; dx + s\; dy, 
q+s\; dx 
+ t\;  dy)=0.\label{rmm}\end{equation} We Taylor expand $\psi$ and use $\psi (x,y,z,p,q)=0$; then 
we see that (\ref{rmm}) is equivalent to
\begin{equation*}\begin{split}\frac{\partial \psi}{\partial x}\cdot dx 
+\frac{\partial 
\psi}{\partial y
}\cdot dy &+\frac{\partial \psi}{\partial z}\cdot (p\; dx +q\; dy 
)\\&+
\frac{\partial \psi}{\partial p}\cdot (r\; dx + s\; dy)
 + \frac{\partial \psi}{\partial q}\cdot  (s\; dx 
+ t\; dy)=0\end{split}\end{equation*}
where the partial derivatives are to be evaluated at $(x,y,z,p,q)$. 
Reorganizing, we see that this is a linear equation system (two 
equations in the three unknowns), and that its augmented matrix is 
the one in (\ref{surff4}).

\medskip

Since at least one of $\frac{\partial \psi}{\partial p}$ and 
$\frac{\partial \psi}{\partial q}$ is invertible, we see that
the rank of the matrix to the left of the augmentation bar in 
(\ref{surff4}) is 2, 
whence it represents a surjective linear map $R^{3}\to R^{2}$;  
the solution set of the equation system is therefore a 1-dimensional 
(and affine) subspace of the $(r,s,t)$-space. So also for a general 
(sufficiently non-singular) ``abstract'' PDE $\Psi \subseteq 
S_{1}(M)$, there are $\infty ^{1}$ solution calottes extending a 
given $P\in \Psi$. 

\medskip

Given $P=(x,y,z,p,q)\in \Psi$. The condition on a neighbour $P'$ that 
it is contained in a calotte $K$ is given by a condition on the 
$(r,s,t)$ of the calotte, namely that it is a solution of the 
equation system (\ref{surff2}) in Proposition \ref{PPrime}; the condition that a 
calotte through $P$ is a solution calotte is that $(r,s,t)$ is a 
solution of the equation system (\ref{surff4}) in Proposition \ref{Psii}. To say 
that $P'$ is a characteristic neighbour of $P$ is therefore to say 
that whenever $(r,s,t)$ solves (\ref{surff4}), it also solves 
(\ref{surff2}). 
 From the ``elementary linear algebra'' in the Appendix therefore 
follows that this is the case iff the 
augmented matrix in (\ref{surff2}) is a scalar multiple of the one in 
(\ref{surff4}).

Therefore we have
\begin{thm}Assume $P=(x,y,z,p,q)$ is in $\Psi$. For $P'= (x+dx, y+dy, 
z+ p\; dx + q\; dy, p+dp,q+dq)$ to be a characteristic neighbour, 
it is necessary and sufficient that there exists a scalar $\lambda$ 
such that 
\begin{equation}(dx,dy,dp,dq) = \lambda \cdot (\psi _{p},\psi _{q}, 
-\psi _{x} - p\cdot \psi _{z} ,-\psi _{y} -q\cdot \psi _{z}  
)\label{reco}\end{equation}
or equivalently, that there exists a scalar $\lambda$ such that
\begin{equation}(dx,dy,dz, dp,dq) = \lambda \cdot (\psi _{p},\psi 
_{q}, p\cdot \psi_{p} + q\cdot \psi_{q},
-\psi _{x} - p\cdot \psi _{z} ,-\psi _{y} -q\cdot \psi _{z}  
)\label{reco2}\end{equation}
\end{thm}
Here, $\psi _{p}$ denotes $\partial \psi /\partial p$ evaluated at 
$P=(x,y,z,p,q)$, and simlarly for $\psi _{q}$, $\psi _{x}$ etc. Note that our assumption that at least one of $\psi _{p}$ and $\psi 
_{q}$ is invertible implies that the scalar $\lambda$ is uniquely 
determined.

\medskip

From the Theorem follows in particular that for $(x+dx,y+dy, z+p\; dx 
+q\; dy)$ to be a characteristic neighbour point of 
$P$ (in the calotte sense), it is necessary that 
\begin{equation}(dx,dy)=\lambda \cdot (\psi _{p},\psi 
_{q})\label{PPP}\end{equation} In fact, it is also sufficient, since the relevant $dp$ and 
$dq$ then can be reconstructed from $\lambda$ and the partial 
derivatives of $\psi$, using (\ref{reco}). 

\medskip

Now that we have an analytic criterion (\ref{reco}) for $P$ and $P'$ 
to be characteristic neighbours, we can also prove that this 
relationship is a symmetric one. The proof is much in the spirit of 
the proof that the relation $\approx$ (``united position'')  is 
symmetric, namely ``ignoring infinitesimals of higher order'':

To say that $P=(x,y,z,p,q)$ and $P' =(x+dx, y+dy,z+dz,p+dp,q+dq)=P+dP$
satisfy $P\approx_{\Psi}P'$ is, by (\ref{reco2}) equivalent to saying 
that
\begin{equation}\label{ag2}dP = \lambda (P;dP)\cdot 
\tilde{\psi}(P),\end{equation}
where we in the $\lambda$-factor record the dependence of the scalar 
$\lambda$ on $P$ as well as on $dP$, and where $\tilde{\psi}:R^{5}\to 
R^{5}$ is the function in the parenthesis on the right hand side of 
(\ref{reco2}), but where we now explicitly record the $P=(x,y,z,p,q)$ 
where the various partial derivatives $\psi _{p}$ etc.\ are to be 
taken. Similarly, to say that $P'\approx _{\Psi}P$ is
 equivalent to saying that
\begin{equation}-dP = \lambda (P';-dP)\cdot 
\tilde{\psi}(P').\label{ag4}\end{equation}
We note that for fixed $P$, and for $dP=0$, we have $\lambda 
(P;dP)=0$. By KL, the function $\lambda (P; -):D(5)\to R$ extends 
uniquely to a linear function $R^{5}\to R$. Since in  the 
expressions on the right hand side of (\ref{ag4}), we have that $dP$ 
occurs linearily, it follows by the ``Taylor principle'' (cf.\ 
\cite{SGM} p.\ 19) -- essentially just Taylor expansion in the direction of 
$dP$ -- that we may replace $P' =P+dP$ by $P$ in both the occurrences 
of $P'$, so that (\ref{ag4}) may be written
$$-dP= \lambda (P; -dP)\cdot \tilde{\psi}(P),$$
which is equivalent to (\ref{ag2}) in view of the linearity of 
$\lambda (P;-)$. This proves the symmetry of $\approx _{\Psi}$.

\medskip
So $\approx _{\Psi}$ is a reflexive symmetric relation on $\Psi$, 
refining the neighbour relation $\sim$. Furthermore, the set of 
$\approx _{\Psi}$-neighbours of a given $P$ in $\Psi$ is 1-dimensional, in a 
sense that  its elements are parametrized by   scalars $\lambda$, 
as is seen in (\ref{reco2}). So $\approx _{\Psi}$ is a 1-dimensional 
geometric distribution, in the sense of \cite{SGM}, \S 2.6.

\subsection*{Differential equation for Monge characteristics}
We consider the surface element $P=(x,y,z,p,q)$ in $\Psi$, so
$\psi (x,y,z,p,q)=0$. A neighbour surface element with same base point 
is of the form $(x,y,z,p+\delta p, q+\delta q)$ with $(\delta p, \delta 
q)\in D(2)$, and this element is in $\Psi$ if $\psi (x,y,z,p+\delta p, 
q+\delta q)=0$; by Taylor expansion, and using $\psi (x,y,z,p,q)=0$, this is 
equivalent to \begin{equation}\label{MMM}(\partial \psi /\partial p )\cdot \delta p + (\partial 
\psi /\partial q)\cdot \delta q =0, \end{equation} where the partial derivatives 
are to be evaluated at $(x,y,z,p,q)$.

A point in $P$ is of the form $(x+dx,y+dy,z+p\; dx + q\; dy)$ with 
$(dx,dy)\in D(2)$, and this point is in the surface  element
$(x,y,z,p+\delta p, q+\delta q)$ iff
$p\; dx + q\; dy = (p+\delta p)\cdot dx + (q+\delta 
q)\cdot dy$, that is, iff\footnote{The reason we did not write $dp$ 
and $dq$, rather than $\delta p$ and $\delta q$ is that notation 
$(dp,dq)$ might lead one  to think that e.g.\  $dx \cdot dp=0$, which we have not 
assumed ; $dx$ and $\delta p$ are what Lie would call {\em 
independent} infinitesimals: $dx \cdot \delta p$ is not assumed to be 
0.}
\begin{equation} dx \cdot \delta p +   dy\cdot \delta q =0.
\label{NNN}\end{equation}
So to say that $(x+dx, y+dy)$ is a Monge-characteristic neighbour of 
$P$ is to say that  all $(\delta p ,\delta q)\in D(2)$ which 
satisfy (\ref{MMM}) also satisfy (\ref{NNN}). Assuming, as before, 
that $\partial \psi /\partial p$ and $\partial \psi /\partial q$ do 
not vanish simultaneously, this property is equivalent to: $(dx,dy)$ 
is of the form $\lambda \cdot (\partial \psi /\partial p, 
\partial \psi /\partial q)$, see Remark after Proposition \ref{abc}. Thus, $(x+dx,y+dy,z+p\; 
dx + q\; dy)$ is a Monge-characteristic neighbour of $P=(x,y,z,p,q)$ 
iff $$(dx,dy)=\lambda \cdot (\partial \psi /\partial p, \partial 
\psi/\partial q).$$
We see that this is just the equation (\ref{PPP}), which is the 
equation for characteristic neighbour point in the calotte sense. We 
conclude that the two notions of ``characteristic neighbour point'' 
agree.

This also proves Theorem \ref{monge}, in a coordinatized situation, 
and since the statement of the theorem is coordinate free, it holds 
also in general.

\section*{Appendix}
Since the linear algebra in question is over the commutative ring $R$ 
which is not a field, only a local ring, we need to elaborate a 
little on the linear algebra/matrix theory over $R$. We use ``vector 
space'' and ``linear'' as synonyms for ``$R$-module'' and 
``$R$-linear. 

Any linear map $R\to R$ is multiplication by a unique $\lambda \in R$.
From this follows, for any vector space $A$:
\begin{prop}\label{abc}Let $p:A\to R$ be a surjective linear map, and 
let $q:A\to R$ 
be any linear map. If the kernel of $p$ is contained in the kernel of 
$q$, then $q=\lambda \cdot p$ for a unique $\lambda \in R$.
\end{prop}
{\bf Proof.} Contemplate the commutative diagram with exact rows
$$\begin{diagram}0&\rTo &Ker (p) &\rTo & A&\rTo ^{p}&R&\rTo 0\\
&&\dTo^{incl}&&\dTo _{id}&&\dTo &\\
0&\rTo &Ker (q) &\rTo & A&\rTo _{q}&R&
\end{diagram}.$$
The right hand vertical map exists by exactness of top row, and is 
multiplication by a unique scalar.

\medskip
\noindent {\bf Remark.} An immediate consequence is that if 
$\va =(a_{1}, \ldots ,a_{n})\in R^{n}$ is a proper vector (meaning: 
at least one of the coordinates 
$a_{i}$ invertible), and if $\vb \in R^{n}$ is a vector such that 
$\va \bullet \vd =0$ implies $\vb \bullet \vd =0$ for all $\vd \in R^{n}$
(where $\bullet$ is the standard dot product of coordinate vectors), 
then $\vb = \lambda \cdot \va$ for some $\lambda \in R$. For this 
conclusion, it  even suffices 
that the implication
$$\va \bullet \vd =0\Rightarrow \vb \bullet \vd =0$$ holds for all 
$\vd \in D(n)$; for, by the KL axiom, a linear map $R^{n}\to R$ is 
completely determined by its value on $D(n)$.
\medskip

Let $p:A\to R$ be a linear map, and let $r\in R$. If $p(x_{0})=r$, then the solution 
set of the equation $p(x)=r$ is the coset $x_{0}+Ker (p)$. As a 
Corollary of the above Propostion, we then have

\begin{prop}Let $p:A\to R$ be a surjective linear map, and $q:A\to R$ 
any linear map. Let $r, s\in R$. If the solution set of $p(x)=r$ is 
contained in the solution set of $q(x)=s$, then there is a unique 
$\lambda \in R$ so that $q=\lambda \cdot p$ and $s=\lambda \cdot r$.
\label{fiv}\end{prop} {\bf Proof.} Take some $x_{0}\in A$ such that 
$p(x_{0})=r$, using $p$ surjective. The assumption then implies that 
$q(x_{0})=s$. The solution sets of the two equations are, 
respectively,
$x_{0}+\ker (p)$, and $x_{0}+\ker (q)$, and the assumed inclusion 
relation then clearly implies $\ker (p) \subseteq \ker (q)$. By the 
previous Proposition, there is a unique $\lambda \in R$ with 
$q=\lambda \cdot p$. We then have
$$\lambda \cdot r = \lambda \cdot p(x_{0}) = q(x_{0})= s.$$

\begin{prop}Consider two linear equation systems
given by the augmented matrices
\begin{equation}\left [ \begin{array}{ccc|c}p_{1}&p_{2}&&r_{1}\\
&p_{1}&p_{2}&r_{2}
\end{array}\right]
\label{first}\end{equation}
and
\begin{equation}\left [ \begin{array}{ccc|c}
q_{1}&q_{2}&&s_{1}\\
&q_{1}&q_{2}&s_{2}
\end{array}\right]
\label{second}\end{equation}
respectively, and assume that at least one of the $p_{i}$s is 
invertible. Assume that the solution set of the first is contained in 
the 
solution set of the second. Then there exists a unique $\lambda \in R$ with
$$\lambda \cdot (q_{1}, q_{2}, s_{1},s_{2}) = 
(p_{1},p_{2},r_{1},r_{2}).$$
\end{prop}
{\bf Proof.} 
Without loss of generality, we may assume that $p_{2}$ is invertible.
Assume $(x,y)$ solves $p_{1}x+p_{2}y= r_{1}$; then there is a 
(unique) $z$ so that $(x,y,z)$  solves the system (\ref{first}). 
Hence by assumption, it solves (\ref{second}), and so $(x,y)$ solves 
$q_{1}x + q_{2}y = s_{1}$. From Proposition \ref{fiv} then follows 
that there exists a $\lambda$ such that 
$q_{1}=\lambda p_{1}$, $q_{2}=\lambda p_{2}$ and $s_{1}=\lambda 
r_{1}$. We prove that also $s_{2}=\lambda r_{2}$: with the unique $z$ 
already considered, we have $p_{1}y+p_{2}z = r_{2}$; multiplying this 
equation  by $\lambda$, we get $q_{1}y + q_{2}z = 
\lambda r_{2}$, but the left hand side here is $s_{2}$ since 
$(x,y,z)$ solves (\ref{second}). This proves the Proposition.

\medskip

If $x$ is a point i  a manifold $M$, the set $\mf _{1}(x)\subseteq M$ 
of first order neighbours of $x$ has a natural ``base'' point, namely 
$x$. This point can be reconstructed from the subset; we claim

\begin{prop}The point $x\in \mf _{1}(x)$ is the only point $z$ with 
the property that for all $y\in M$ with $y\sim z$, we have $y\in \mf 
_{1}(x)$.
\label{basept}\end{prop}
{\bf Proof.} Since the assertion is coordinate free, it suffices to 
prove it for the case where $M=R^{n}$ and $x=0 \in R^{n}$. Note that 
now $\mf _{1}(x) =D(n)$. Then the assertion of the Proposition 
amounts to the assertion: if $z\in D(n)$ has the property that $z+u 
\in D(n)$ for all $u\in D(n)$, then $z=0$.
To prove that the first coordinate $z_{1}$ of $z$ is $0\in R$, we 
use that
$z+(d,0,\ldots ,0) \in D(n)$ for all $d\in D$, which implies
$(z_{1}+d)^{2}=0$ for all $d\in D$.
Now
$$0=(z_{1}+d)^{2} = z_{1}^{2}+d^{2}+2z_{1}d = 2z_{1}d.$$
Since this holds for all $d\in D$, it follows from KL that 
$2z_{1}=0$, hence $z_{1}=0$. Similarly for the other coordinates 
$z_{2}, \ldots ,z_{n}$.

\medskip

There are similar characterizations of $\mf_{k}(x)$ as as a subset of 
$\mf _{k+1}(x)\subseteq M$ for $k=2, \ldots $. One then needs that the 
integer $k+2$ is invertible in $R$.

\medskip

\small

\noindent University of Aarhus, November 2010

\noindent
\begin{verbatim}kock@imf.au.dk\end{verbatim}

\end{document}